\documentclass{elsart3-1}

\usepackage{amssymb,amsfonts,verbatim}
\usepackage[tbtags]{amsmath}

\usepackage[english,francais]{babel}

\newtheorem{theorem}{Theorem}[section]

\newtheorem{e-proposition}[theorem]{Proposition}

\newtheorem{e-definition}[theorem]{Definition\rm}

\renewcommand{\theequation}{\arabic{equation}}
\def\og{\leavevmode\raise.3ex\hbox{$\scriptscriptstyle\langle\!\langle$~}}
\def\fg{\leavevmode\raise.3ex\hbox{~$\!\scriptscriptstyle\,\rangle\!\rangle$}}

\journal{}
\begin{document}
\centerline{}
\begin{frontmatter}


\selectlanguage{english}
\title{Remarks on local controllability for the Boussinesq system with Navier boundary condition 
	}


\selectlanguage{english}
\author[authorlabel1]{Cristhian Montoya}
\ead{cristhian.montoya@usm.cl}

\address[authorlabel1]{Universidad T\'ecnica Federico Santa Maria, Casilla 110--V, Valparaiso, Chile}


\medskip

\begin{abstract}
\selectlanguage{english}
	This note deals with the local exact controllability to a particular class of trajectories
    for the Boussinesq system with nonlinear Navier--slip boundary conditions and internal controls 
    having vanishing components.  Briefly speaking, in two dimensions, the local exact controllability 
    property is obtained using only one control in the heat equation,  meanwhile  two scalar controls are required in 
    three dimensions.

\vskip 0.5\baselineskip

\selectlanguage{french}
\noindent{\bf R\'esum\'e} \vskip 0.5\baselineskip \noindent
{\bf Remarque sur la contr\^{o}labilit\'e locale du syst\`{e}me de Boussinesq avec la condition de fronti\`{e}re de Navier.} 
	Cette note concerne la contr\^{o}labilit\'e locale d'une classe particuli\`{e}re de trajectoires, ceci pour le syst\`{e}me 
	de Boussinesq avec la condition de Navier non lin\'{e}aire et certains contr\^{o}les internes. Bri\`{e}vement, la 
	propi\'et\'e de contr\^{o}labilit\'e exacte locale s'obtient en dimension deux n'utilisant que le contr\^{o}l associ\'e
	{\`{a}} l'\'{e}quation de la chaleur. Tandis que, deux contr\^{o}les scalaires sont n\'ecessaires pour obtenir n\^{o}tre 
	r\'esultat dans le cas de dimension trois.


\end{abstract}
\end{frontmatter}

\selectlanguage{french}

\selectlanguage{english}
\section{Introduction}\label{section.intro}
	The interaction of incompressible fluids with  a diffusion process can be modeled by a coupled 
	system between the Navier--Stokes and heat equations, usually called Boussinesq system. 
	On bounded domains, both heat and the velocity field can show a different behaviour on its 
	boundary. In this paper, nonlinear Navier--type boundary conditions for the fluid flow and homogeneous Neumann
	conditions for the diffusion equation are considered in order to study the local exact controllability 
	for the Boussinesq system with few scalar controls.
		
	Henceforth, let $\Omega$  be a nonempty bounded connected open subset of $\mathbb{R}^N$ ($N=2$ or $N=3$) of class 
    $C^{\infty}$. Let $T > 0$ and let   $\omega\subset \Omega$ be a (small) nonempty open subset which is    
    the control domain. Here, we will use the notation $Q:=\Omega\times(0,T)$,\, $\Sigma:=\partial\Omega\times(0,T)$ 
    and  $n$ the outward unit normal vector to $\Omega$. Moreover, $C$ denotes a generic positive constant
    which may depend on $\Omega$ and $\omega$.
    
    In this Note, we will consider the Boussinesq system with Navier--slip  and Neumann conditions
	\begin{equation}\label{intro_mainsystem}
    \left\{
    \begin{array}{llll}
        y_{t}-\nabla\cdot (Dy)+(y,\nabla)y+\nabla p=u\chi_{\omega}+\theta e_{N},\quad  \nabla \cdot y=0  &\text{ in }& Q,\\
        \theta_t-\Delta\theta+y\cdot\nabla\theta=v1_{\omega}&\text{ in }& Q,\\
        y\cdot n=0,\, (\sigma(y,p)\cdot n)_{tg}+f(y)_{tg}=0,\quad  \nabla\theta\cdot n=0 &\text{ on }&\Sigma,\\
        y(\cdot,0)=y_0(\cdot),\,\,\theta(\cdot,0)=\theta_0(\cdot) & \text{ in }&\Omega,
    \end{array}
    \right.
	\end{equation}
	as well as the linearized Boussinesq system (around a target flow of the form 
	($(0,\overline{p},\overline{\theta})$)
    \begin{equation}\label{intro_linearsys}
    \left\{
    \begin{array}{llll}
        y_{t}-\nabla\cdot (Dy)+\nabla p=h_1+u\chi_{\omega}+\theta e_{N},\quad  \nabla \cdot y=0  &\text{ in }& Q,\\
        \theta_t-\Delta\theta+y\cdot\nabla\overline{\theta}=h_2+v1_{\omega}&\text{ in }& Q,\\
        y\cdot n=0,\, (\sigma(y,p)\cdot n)_{tg}+(A(x,t)y)_{tg}=0,\quad \nabla\theta\cdot n=0&\text{ on }&\Sigma, \\
        y(\cdot,0)=y_0(\cdot),\,\,\theta(\cdot,0)=\theta_0(\cdot) & \text{ in }&\Omega,
    \end{array}
    \right.
	\end{equation}       
     where $y=y(x,t)$ is the velocity field of the fluid, $\theta=\theta(x,t)$ their temperature, 
     $v$ and $u=(u_1,\dots, u_N)$ stands for the controls, which are acting in a arbitrary fixed domain $\omega\times(0,T)$,
     where $\chi_{\omega}$ is a smooth positive function such that $\chi_{\omega}=1$ in $\omega'$, 
     $\omega'\Subset\omega$, and $1_{\omega}$ is the indicator function. Here, the gravity vector field is given by
     $e_N=(0,1)$ for $N=2$, or $e_N=(0,0,1)$ for $N=3$. Moreover, $f:\mathbf{R}^{N}\to \mathbf{R}^{N}$ is a 
     nonlinear regular function given, $\sigma(y,p):=-pId+Dy$ is the stress tensor, $A$ is a $N\times N$ 
     matrix--valued function in a suitable space, and {\textit{tg}} stands for the 
     tangential component of the corresponding vector field, i.e.,  $y_{tg}=y-(y\cdot n)n$.
          
    In the context of controllability, the first results for the Boussinesq system 
    were made by Fursikov and Imanuvilov in \cite{fursikov1998local} and \cite{fursikov1999exact}. 
    The work by S. Guerrero \cite{guerrero2006local_boussinesq} shows the local exact controllability to the trajectories 
    of the Boussinesq system with Dirichlet boundary conditions and  $N+1$ distributed scalar controls supported in small sets.
    
    Additionally, recent works have been developed for controllability problems with reduced number of controls. 
    For instance,  N. Carre\~no and S. Guerrero in \cite{carreno2013local} have proven the local null controllability for the 
    Navier--Stokes with Dirichlet conditions and $N-1$ scalar controls. The recent work made by S. Guerrero and 
    C. Montoya shows that the local null controllability property is achieved for the $N$--dimensional Navier--Stokes system
    with Navier--slip conditions and $N-1$ scalar controls \cite{guerrero2018local}. The methodology in the previous articles
    are Carleman estimates. 
    In the three dimensional case of the Navier--Stokes system with Dirichlet conditions, J-M. Coron and P. Lissy developed
    in \cite{coron2014local} a new strategy to prove the local null controllability using only one scalar control. 
    
    Concerning the  $N$-dimensional Boussinesq system with Dirichlet conditions,  in  \cite{fernandez2006some} 
    the authors proved that the local exact controllability
    to the trajectories can be achieved with $N-1$ scalar controls, under certain geometric assumption on the control 
    domain. N. Carre\~no showed the 
    local controllability of the $N$--Boussinesq system using $N-1$ scalar controls, without conditions on the control 
    domain \cite{carreno2012local}. Finally, this Note improves the results of \cite{carreno2013local} and 
    \cite{guerrero2018local}.

	Our results below extend the results of \cite{carreno2013local} and \cite{guerrero2018local}. Taking into account the 
	relation between the observability and controllability property, it will be appropriate to consider the following adjoint 
	system related to \eqref{intro_linearsys}:
	\begin{equation}\label{3.adjointsystem}
    \left\{
    \begin{array}{llll}
        -\varphi_{t}-\nabla\cdot (D\varphi)+\nabla \pi=g-\psi\nabla\overline{\theta},\quad  \nabla \cdot \varphi=0
        &\text{ in }& Q,\\
        -\psi_t-\Delta\psi=g_0+\varphi\cdot e_{N}&\text{ in }& Q,\\
        \varphi\cdot n=0,\, (\sigma(\varphi,\pi)\cdot n)_{tg}+(A^t(x,t)\varphi)_{tg}=0,\quad \nabla\psi\cdot n =0
        &\text{ on }&\Sigma, \\
        \varphi(\cdot,T)=\varphi^T(\cdot),\,\,\psi(\cdot,T)=\psi^T(\cdot) & \text{ in }&\Omega,
    \end{array}
    \right.
	\end{equation} 
	 where $g, \varphi^T, g_0$ and $\psi^T$ satisfying adequate regularity assumptions.   
	We will introduce several spaces and hypotheses over $\overline{\theta}$ which will be needed in order to have 
	suitable Carleman estimates for the solution of
    \eqref{3.adjointsystem}:
    \begin{equation*}
    \begin{array}{ll}
    	W&=\{u\in H^1(\Omega)^N: \nabla\cdot u=0 \,\,\text{in}\,\, \Omega,\,\,\, u\cdot n=0
     	\,\,\text{on}\,\, \partial\Omega\},\\
     	H&=\{ u\in L^{2}(\Omega)^{N}: \nabla\cdot u=0, \,\,\text{in}\,\, \Omega\,\,\, u\cdot n=0 \,\,\text{on}\,\,
      	\partial\Omega \},\\
      	 P^1_{\varepsilon}&=H^{5/4+\varepsilon}{(0,T;L^2(\partial\Omega)^{N\times N}}),\quad
    	P^{2}=L^2(0,T;H^{5/2}(\partial\Omega)^{N\times N}),\quad \forall \varepsilon >0,\\
    	Y_m&:=L^2(0,T;H^{2m}(\Omega)^N)\cap H^m(0,T;L^2(\Omega)^N),\quad m=1,2. 	
    \end{array}
	\end{equation*}
	and 
	\begin{equation}\label{intro_reg_trajectheta}
    \overline{\theta}\in L^{\infty}(0,T;W^{3,\infty}(\Omega)),\quad
    \nabla\overline{\theta}_t \in L^{\infty}(Q)^N.
	\end{equation}
    Here, the target flow $(0,\overline{p},\overline{\theta})$ satisfies the problem
	\begin{equation}\label{intro_trajectorysys2}
    \left\{
    \begin{array}{llll}
        \nabla \overline{p}=\overline{\theta} e_{N},\quad \overline{\theta}_t-\Delta\overline{\theta}=0  &\text{ in }& Q,\\
        \nabla\overline\theta\cdot n =0&\text{ on }&\Sigma,\\
        \overline{\theta}(\cdot,0)=\overline{\theta}_0(\cdot) & \text{ in }&\Omega.
    \end{array}
    \right.
	\end{equation}
	Our first main result is a new Carleman estimate for the solution of \eqref{3.adjointsystem}. Several weight functions 
	are needed:
	\begin{equation}\label{3.weights}
    \begin{array}{llllll}
    \alpha(x,t) &= \dfrac{e^{2\lambda\|\eta\|_{\infty}}-e^{\lambda\eta(x)}}
    {(t(T-t))^{11}}, \quad & \xi(x,t)&=\dfrac{e^{\lambda\eta(x)}}{(t(T-t))^{11}},\quad &
    \alpha^{*}(t) &= \max_{x\in\overline{\Omega}} \alpha(x,t),\\
    \xi^{*}(t) &= \min_{x\in\overline{\Omega}} \xi(x,t),\quad &
    \widehat\alpha(t) &=\min_{x\in\overline{\Omega}} \alpha(x,t),\quad &
    \widehat\xi(t) &= \max_{x\in\overline{\Omega}} \xi(x,t).
    \end{array}
	\end{equation}
	Here, $\eta\in C^2(\overline{\Omega})$ and satisfies that    
    \[|\nabla \eta|>0 \mbox{ in }\overline{\Omega}\setminus\omega_0,\,\,\,\, \eta>0
     \mbox{ in }\Omega\,\,\, \mbox{ and }\,\, \eta \equiv 0 \mbox{ on }\partial
     \Omega,\]
     where  $\omega_0 \Subset \omega_1\Subset\omega'\Subset{\omega}$ is a nonempty open set. The existence of such a 
     function $\eta$ is proved in \cite{fursikov1996controllability}.
     \vskip 0.2cm
 
	\begin{theorem}\label{th.carleman}
 		Assume $A\in P^1_{\varepsilon}\cap P^2$ and $(0,\overline{p},\overline\theta)$ satisfying 
 		\eqref{intro_reg_trajectheta}--\eqref{intro_trajectorysys2}. There exists a constant $\lambda_0$, such that for any 
   		$\lambda\geq\lambda_0$  there exist two constants $C(\lambda)>0$ increasing on 
    	$\|A\|_ {P^1_{\varepsilon}\cap P^2}$ and $s_0(\lambda)>0$ such that for any $j\in\{1,2\}$, any $a>0$, any 
    	$g\in L^2(Q)^3$, any $g_0\in L^2(Q)$, any $\varphi^T\in H$ and any $\psi^T\in L^2(\Omega)$, the solution
   		of \eqref{3.adjointsystem} satisfies
    \begin{equation}\label{3. carleman_estimateN3}
        \begin{aligned}
        s^3&\iint\limits_{Q}e^{-2(1+a)s\alpha^*}(\xi^*)^3|\varphi|^2dxdt
        +s^5\iint\limits_{Q}e^{-2(1+a)s\alpha^*}(\xi^*)^5|\psi|^2dxdt\\
        &\leq C\Bigl(\iint\limits_{Q}e^{-2as\alpha^*}(|g|^2+|g_0|^2)dxdt
        +(N-2)s^7\int\limits_{0}^{T}\int\limits_{\omega'}e^{-4s\hat{\alpha}+2(1-a)s\alpha^*}
        (\hat{\xi})^{12}|\varphi_j|^2dxdt\\
        &\hspace{2cm}+s^{13}\int\limits_{0}^T\int\limits_{\omega}
        e^{-8s\hat\alpha+(6-2a)s\alpha^*}(\hat\xi)^{24}|\psi|^2 dxdt 
        \Bigr)
        \end{aligned}
    \end{equation}
    for every $s\geq s_0$.  
    \end{theorem}
  
	The second main result in this Note concerns the local controllability to a particular class of trajectories 
 	of \eqref{intro_mainsystem}. This result is presented as follows:
 	\vskip 0.2cm

 	\begin{theorem}\label{maintheorem}
 		Assume $f\in C^4(\mathbb{R}^N;\mathbb{R}^N)$ with $f(0)=0$ and $i\in\{1,\dots, N-1\}$ fixed. Let   
    	$(0,\overline{p},\overline{\theta})$ be a solution to \eqref{intro_trajectorysys2} satisfying 
    	\eqref{intro_reg_trajectheta}. Then, for every $T>0$ and $\omega\subset\Omega$, there exists $\delta>0$ 
    	such that, for every $(y_0,\theta_0)\in [H^3(\Omega)^N\cap W]\times H^1(\Omega)$ satisfying 
    	\begin{equation}\label{intro_compatib_cond}
        	(Dy_0\cdot n)_{tg}+(f(y_0))_{tg}=0\mbox{ on }\partial\Omega\quad\mbox{and}\quad
        	|(y_0,\theta_0)-(0,\overline{\theta}_0)\|_{[H^3(\Omega)^N\cap W]\times H^1(\Omega)}\leq \delta, 
    	\end{equation}
   		 we can find  controls $v\in L^2(\omega\times(0,T))$ and 
    	$u\in L^2(0,T;H^2(\omega)^N)\cap H^1(0,T;L^2(\omega)^N)$ with $u_i\equiv 0$ and $u_N\equiv 0$ 
    	such that the corresponding solution $(y,p,\theta)$ to \eqref{intro_mainsystem} satisfies
		\begin{equation}\label{intro_def_localexactc}
    	y(\cdot,T)=\overline{y}(\cdot,T)\quad\mbox{and}\quad \theta(\cdot,T)=\overline{\theta}(\cdot,T)
    	\quad\mbox{in}\,\,\Omega.
		\end{equation}
 	\end{theorem}
 	In the following sections, we will indicate the main ideas of the proof of Theorem \ref{th.carleman} and 
 	Theorem \ref{maintheorem}. 
\section{\normalsize A new Carleman inequality}\label{section2}
	In this section, we give the proof  of Theorem \ref{th.carleman}. Our arguments are based in 
	\cite{carreno2013local,CoronGuerrero09,fernandez2006some,guerrero2018local}. 
	From \eqref{3.adjointsystem} and using the decomposition $\rho\varphi=w+z,\quad \rho\pi=\pi_{z}+\pi_{w}$ and 
	$\rho\psi=\tilde{\psi}$, where $\rho(t)=e^{-as\alpha^*}$ and $a>0$, it is very easy to verify that 
	$(w,\pi_{w}), (z,\pi_z)$ and $\tilde{\psi}$ are solutions to the systems
  \begin{equation*}\label{3.n3_system_w}
    \left\{
    \begin{array}{lllll}
        -w_{t}-\nabla\cdot (Dw)+\nabla \pi_{w}=\rho g; &\quad 
        -z_{t}-\nabla\cdot (Dz)+\nabla \pi_{z}=-\rho' \varphi-\tilde{\psi}\nabla\overline{\theta} &\text{ in }& Q,\\
        \nabla \cdot w=0; &\quad \nabla \cdot z=0 &\text{ in }& Q,\\
        w\cdot n=0,\, (\sigma(w,\pi_{w})\cdot n)_{tg}+(A^t(x,t)w)_{tg}=0; &\quad 
         z\cdot n=0,\, (\sigma(z,\pi_{z})\cdot n)_{tg}+(A^t(x,t)z)_{tg}=0 &\text{ on }&\Sigma, \\
        w(\cdot,T)=0;&\quad   z(\cdot,T)=0& \text{ in }&\Omega,
    \end{array}
    \right.
	\end{equation*}  
    and 
	\begin{equation}\label{3.n3_system_tildepsi}
    \left\{
    \begin{array}{llll}
        -\tilde{\psi}_{t}-\Delta\tilde{\psi}=\rho g_0+\rho\varphi_3-\rho'\psi &\text{ in }& Q,\\
        \nabla\tilde{\psi}\cdot n=0 &\text{ on }&\Sigma, \\
        \tilde{\psi}(\cdot,T)=0 & \text{ in }&\Omega,
    \end{array}
    \right.
	\end{equation}  
	We will use the Carleman inequality for parabolic equations with Neumann conditions \cite{fursikov1996controllability} 
	for the system  \eqref{3.n3_system_tildepsi} in  order to estimate the global terms associated to $\tilde\psi$. 
    Thus, there exists $\tilde{\lambda}>0$ such that for any $\lambda>\tilde{\lambda}$ there exists  a positive constant $C$ 
    depending on  $\lambda,\Omega,\omega_2,\|\overline{\theta}\|_{L^{\infty}(0,T;W^{3,\infty}(\Omega))}$ such that
	\begin{equation}\label{3.n3_carleman_linear_Fourier}
    \begin{aligned}
       \iint\limits_{Q}&e^{-2s\alpha}(s\xi|\tilde{\psi}_t|^2
       +s\xi\sum\limits_{\ell,m=1}^3|\partial_{\ell m}\tilde{\psi}|^2+ s^3\xi^3
       |\nabla \tilde{\psi}|^2+s^5\xi^5|\tilde{\psi}|^2)dxdt\\
        \leq & C\Bigl(\iint\limits_{Q}e^{-2s\alpha}s^2\xi^2(|\rho g_0|^2+|\varphi_3|^2
        +|\rho'|^2|\rho|^{-2}|\tilde{\psi}|^2) dxdt  
        +s^5\int\limits_{0}^T\int\limits_{\omega_1}e^{-2s\alpha}\xi^5|\tilde{\psi}|^2 dxdt\Bigr),
    \end{aligned}
	\end{equation}
    for every  $s\geq C$.
    
    The arguments below are given for the case $N=3$. For $k=1,3$, we can deduce the inequality
   	\begin{equation}\label{3.n3_end_step2}
    \begin{aligned}
        &I(s,z)+J(s,\tilde{\psi})
         \leq C\Biggl( \|\rho g\|^2_{L^2(Q)^3}+\|\rho g_0\|^2_{L^2(Q)}
        +s^5\int\limits_{0}^T\int\limits_{\omega_1}e^{-2s\alpha}\xi^5|\tilde{\psi}|^2 dxdt\\
        &\quad +\sum\limits_{k=1,3}\Biggl[\int\limits_{0}^T\int\limits_{\omega_1}e^{-2s\alpha}(s^5\xi^5|
        z_k|^2+s^3\xi^3|\nabla z_k|^2)dxdt
        +\int\limits_{0}^T\int\limits_{\omega_2}e^{-2s\alpha}\xi^2|\nabla\partial_{k}
        \pi_{z}|^2 dxdt
        \Biggr] \Biggr),    
    \end{aligned}
	\end{equation}
    where $J(s,\tilde{\psi})$ denotes the left--hand side of \eqref{3.n3_carleman_linear_Fourier}, and for $k=1,3$,  $I(s,z)$ is 
    defined by 
    \begin{equation*}\label{3.n3_def_I}
    \begin{aligned}
        I(s, z):=\sum\limits_{k=1,3}&s^5\iint\limits_{Q}e^{-2s\alpha}\xi^5|z_k|^2dx dt
        +s^3\iint\limits_{Q}e^{-2s\alpha}\xi^3|\nabla z_k|^2dx dt
        +s^3\iint\limits_{Q}e^{-2s\alpha}\xi^3|z_2|^2dx dt\\
        &+\|s^{1/2}e^{-s\alpha^*}(\xi^*)^{9/22}z\|^2_{Y_1}
        +\|s^{-1/2}e^{-s\alpha^*}(\xi^*)^{-15/22}z\|^2_{L^2(0,T;H^4(\Omega)^3)\cap H^2(0,T;L^2(\Omega)^3)}\\
        &+\|s^{1/2}e^{-s\alpha^*}(\xi^*)^{9/22}{\pi}_z\|^2_{L^2(0,T;H^1(\Omega))}.
     \end{aligned}
	\end{equation*}
  	Here, $\omega_1$ and $\omega_2$ are open sets such that  $\omega_1\Subset \omega_2\Subset \omega$. The rest of the proof
  	is oriented towards the absorption of the local pressure term in \eqref{3.n3_end_step2}. However, we have omitted these 
  	details since analogous arguments can be found in \cite{guerrero2018local}, Section 3. Let us remark that the regularity 
  	over $\overline{\theta}$ given in \eqref{intro_reg_trajectheta} is used in several estimates associated to the 
  	pressure term. The other local terms can be 
  	estimated in an easier way. Therefore, those local estimates lead to the desired Carleman inequality 
  	\eqref{3. carleman_estimateN3}.
  	
\section{\normalsize Local controllability for the Boussinesq system}\label{section2}
	The proof of Theorem \ref{maintheorem} follows the ideas in \cite{carreno2013local} and \cite{guerrero2018local}. 
	Thus, in a first step a null controllability result for \eqref{intro_linearsys} with an appropriate right--hand side
	$h_1,\, h_2$. Here, the idea is to look for a solution in an appropriate weighted functional space. Let us 
    \[L_1w:=w_t-\nabla\cdot Dw\quad\mbox{and}\quad L_2w:=w_t-\Delta w\]
    and let us define the space $E$ as follows:
	\begin{equation*}
    \begin{aligned}
        &\{(y,p,u_1,\theta, v): e^{as\beta^*}y, 
        e^{2s\hat{\beta}-(1-a)s\beta^*}(\hat{\gamma})^{-6}(u_1,0,0)\chi_{\omega},\,  
        \tilde{\rho}(\partial_t u_1,0,0)\in L^2(Q)^3,e^{as\beta^*}\theta\in L^2(Q),\\
        &\hspace{5mm}e^{4s\hat\beta-(3-a)s\beta^*}(\hat\gamma)^{-12}v1_{\omega}\in L^2(Q),\,
        \tilde{\rho}u_1\in L^2(0,T;H^2(\Omega)),\,supp\,\, u_1\subset\omega\times(0,T),\,\\ 
        &\hspace{5mm}e^{as\beta^*}(\gamma^*)^{-12/11}y\in Y_1,\, 
        e^{as\beta^*}(\gamma^*)^{-12/11}\theta\in L^2(0,T;H^2(\Omega))\cap H^1(0,T;L^2(\Omega)),\\
        &\hspace{5mm}e^{(a+1)s\beta^*}(\gamma^*)^{-3/2}
        (L_1y+\nabla p-(u_1,0,0)\chi_{\omega}-\theta e_3)\in L^2(Q)^3,\,\\
        &\hspace{5mm}e^{(a+1)s\beta^*}(\gamma^*)^{-5/2}
        (L_2\theta+y\cdot\nabla\overline\theta-v1_{\omega})\in L^2(Q)
         \}=:E,
    \end{aligned}
	\end{equation*}
    where $\tilde{\rho}:= e^{4s\hat{\beta}+2(1-a)s\beta^*}(\hat{\gamma})^{-12}e^{-(1+a)s\beta^*}(\gamma^*)^{9/22}$ 
    and whose weight functions are given by 
     \begin{equation}\label{3.3_weight_not_vanishing}
        \begin{aligned}
            &\beta(x,t) = \dfrac{e^{2\lambda\|\eta\|_{\infty}}-e^{\lambda\eta(x)}}
            {\ell^{11}(t)},\quad 
            &\gamma(x,t)=\dfrac{e^{\lambda\eta(x)}}{\ell^{11}(t)},\quad 
            &\beta^*(t) = \max_{x\in\overline{\Omega}} \beta(x,t),\\ 
           &\gamma^*(t) = \min_{x\in\overline{\Omega}} \gamma(x,t),\quad
            &\widehat\beta(t) = \min_{x\in\overline{\Omega}} \beta(x,t),\quad  
            &\widehat\gamma(t) = \max_{x\in\overline{\Omega}} \gamma(x,t).
        \end{aligned}
    \end{equation}
	In this case, $\ell\in C^2([0,T])$ is a positive function in $[0,T)$ such that $\ell(t)> t(T-t)$ for all 
     $t\in [0,T/4]$ and $\ell(t)=t(T-t)$ for all $t\in [T/2,T]$.
 	 \vskip 0.2cm

 	\begin{e-proposition}\label{3.3_pro_controllability}
     Let $s$ and $\lambda$ be like in Theorem \ref{th.carleman} and 
    $(0,\overline{p},\overline{\theta})$ satisfy \eqref{intro_trajectorysys2}. Assume that
     \begin{equation}\label{3.3_condition_data}
        y_0\in W, \theta_0\in H^1(\Omega), \,\, e^{(a+1)s\beta^*}(\gamma^*)^{-3/2}h_1\in L^2(Q)^3
        \,\,\,\mbox{and}\,\, e^{(a+1)s\beta^*}(\gamma^*)^{-5/2}h_2\in L^2(Q).     
    \end{equation}
    Then, there exists controls $u_1$ and $v$ such that, if $(y,p,\theta)$ is the associated solution to 
    \eqref{intro_linearsys}, we have $(y, p, u_1,\theta,v)\in E$. In particular $y(\cdot, T)=0$ and 
    $\theta(\cdot,T)=0$ in $\Omega$. 
    \end{e-proposition}
    \vskip 0.2cm

  	The rest of the proof of Theorem \ref{maintheorem} relies on two fixed point theorems, namely, one for the nonlinearity posed
  	on the boundary condition, and another one, for the convective term in \eqref{intro_mainsystem}. We will
  	mention only these results since the methodology given in \cite{guerrero2018local} can be adapted to  
  	\eqref{intro_mainsystem}. Thus, for $N=3$, we consider the nonlinear system
  	\begin{equation}\label{4.system_nonlinearity_boundary}
    \left\{
    \begin{array}{llll}
         y_{t}-\nabla\cdot (Dy)+\nabla p=h_1+(u_1,0,0)\chi_{\omega}+\theta e_{3},\quad  \nabla \cdot y=0  &\mbox{ in }& Q,\\
        \theta_t-\Delta\theta+y\cdot\nabla\overline{\theta}=h_2+v1_{\omega}&\mbox{ in }& Q,\\
        y\cdot n=0,\, (\sigma(y,p)\cdot n)_{tg}+(f(y))_{tg}=0,\quad  \nabla\theta\cdot n=0&\mbox{ on }&\Sigma, \\
        y(\cdot,0)=y_0(\cdot),\,\,\theta(\cdot,0)=\theta_0(\cdot) & \mbox{ in }&\Omega.
    \end{array}
    \right.
	\end{equation}
	\begin{theorem}\label{4.teo_nonlinear_boundary}  
    Let us assume that $f\in C^4(\mathbf{R}^3;\mathbf{R}^3)$ with 
    $f(0)=0$. Then, for every $T>0$ and $\omega\subset\Omega$, there exists 
    $\delta>0$ such that, for every $a>0$ and for every $(y_0,\theta_0)\in H^3(\Omega)^3\cap W\times H^1(\Omega)$,\, 
    $h_1\in Y_1, h_2\in L^2(Q)$ satisfying  
    $e^{(a+1)s\beta^*}(\gamma^*)^{-3/2}h_1\in L^2(Q)^3$\, and\, $e^{(a+1)s\beta^*}(\gamma^*)^{-5/2}h_2\in L^2(Q),$
    \begin{equation}\label{4.1.small_data} 
    \|h_1\|_{Y_1}+\|h_2\|_{L^2(Q)}+\|y_0\|_{H^3(\Omega)^3\cap W}+\|\theta_0\|_{H^1(\Omega)}\leq \delta
    \end{equation}
    and (\ref{intro_compatib_cond}), there exists controls $v\in L^2(0,T;L^2(\omega))$ and
    $u_1\in L^2(0,T;H^2(\omega))\cap H^1(0,T;L^2(\omega))$ and an associated 
    solution $(y,p,\theta)$ of (\ref{4.system_nonlinearity_boundary}) satisfying 
    $(y,\theta)\in Y_2\times L^2(0,T;H^2(\Omega))\cap H^1(0,T;L^2(\Omega))$ and 
    such that $(y,p,u_1,\theta,v)\in E$.
	\end{theorem}
	\vskip 0.2cm

	\begin{theorem}\label{4.2.teo_inverse_mapping}
    Suppose that $\mathcal{B}_1,\mathcal{B}_2$ are Banach spaces and 
    \[\mathcal{A}:\mathcal{B}_1\to \mathcal{B}_2\] is a continuously differentiable map. We assume that for 
    $b_1^0\in \mathcal{B}_1, b_2^0\in \mathcal{B}_2$ the equality
    \begin{equation}\label{4.2.def_b_2^0}
    \mathcal{A}(b_1^0)=b_2^0
    \end{equation}
    holds and $\mathcal{A}'(b_1^0):\mathcal{B}_1\to \mathcal{B}_2$ is an epimorphism. Then there exists 
    $\delta >0$ such that for any $b_2\in \mathcal{B}_2$ which satisfies the condition
    $$\|b_2^0-b_2\|_{\mathcal{B}_2}<\delta$$
    there exists a solution $b_1\in \mathcal{B}_1$ of the equation $$\mathcal{A}(b_1)=b_2.$$
	\end{theorem}    
    Let us set 
    $$y=\tilde{y},\,\, p=\overline{p}+\tilde{p}\,\,\mbox{and}\,\, \theta=\overline{\theta}+\tilde{\theta}.$$
    For $a=2>1$, we apply Theorem \ref{4.2.teo_inverse_mapping} with the spaces
    $$\mathcal{B}_1:=\{(y,p,u_1,\theta,v)\in E: y\in Y_2\},$$
    $$\mathcal{B}_2:=\{(h_1,y_0,h_2,\theta_0)\in Z_1\times [H^3(\Omega)^3\cap W]\times Z_2\times H^1(\Omega):    
    h_1,h_2,y_0,\theta_0\,\,\mbox{satisfies}\,\,(\ref{4.1.small_data}\}),$$
    and where $$Z_1:=L^2(e^{3s\beta^*}(\gamma^*)^{-3/2}(0,T);L^2(\Omega)^3),\quad\mbox{and}\quad  
    Z_2:=L^2(e^{3s\beta^*}(\gamma^*)^{-5/2}(0,T);L^2(\Omega)).$$
    By defining the operator $\mathcal{A}:\mathcal{B}_1\to \mathcal{B}_2$ by
    $$\mathcal{A}\to (L_1\tilde{y}+(\tilde{y}\cdot\nabla)\tilde{y}
    +\nabla\tilde{p}-\tilde{\theta}e_3-(u_1,0,0)\chi_{\omega},\tilde{y}_0,L_2\tilde{\theta}
    +\tilde{y}\cdot\nabla\tilde{\theta}+\tilde{y}\cdot\nabla\overline{\theta}-v1_{\omega},\tilde{\theta}_0),$$
    for every $(\tilde{y},\tilde{p},u_1,\tilde{\theta},v)\in\mathcal{B}_1$, one can easily check the conditions
    for $\mathcal{A}$ in order to complete the proof of Theorem \ref{maintheorem}. 
    \vskip 0.3cm

 	\noindent\textbf{Some open problems.} It would be interesting to know if the local controllability to 
 	the trajectories with $N-1$ scalar controls holds for $\overline{y}\neq 0$ and $\omega$ like in Theorem \ref{maintheorem}. 
 	However,  is not clear at all and therefore is an open problem even for the Navier--Stokes system.
     
    On the other side, could be reasonable to expect results of the same kind whether one considers nonlinear 
    conditions such as $\nabla\theta\cdot n+g(\theta)=0$, where $g$ is a suitable function to study.
    
    Recently, Coron et al. have proved a global exact controllability result  for 
    the Navier--Stokes and Navier--type conditions (for small time), see  \cite{2018coronmarbach}. 
    A challenging problem would be to use the Boussinesq system proposed in this Note in order to apply and prove analogous
    results to \cite{2018coronmarbach}.

\section*{Acknowledgements} The author would like to express his gratitude to Sergio Guerrero for his suggestions,
	 which have contributed to a better presentation of this paper. 
	 
	This work has been supported by FONDECYT grant 3180100. 


\begin{thebibliography}{00}

 \bibitem{carreno2013local}
N.~Carre{\~n}o and S.~Guerrero.
\newblock Local null controllability of the n-dimensional Navier--Stokes system
  with n- 1 scalar controls in an arbitrary control domain.
\newblock {\em Journal of Mathematical Fluid Mechanics}, 15(1):139--153, 2013.

\bibitem{carreno2012local}
N.~Carre{\~n}o.
\newblock Local controllability of the n-dimensional Boussinesq system with n-1
  scalar controls in an arbitrary control domain.
\newblock {\em arXiv preprint arXiv:1201.1871}, 2012.

\bibitem{CoronGuerrero09}
J.-M. Coron and S.~Guerrero.
\newblock Null controllability of the n-dimensional Stokes system with n- 1
  scalar controls.
\newblock {\em Journal of Differential Equations}, 246(7):2908--2921, 2009.

\bibitem{coron2014local}
J.-M. Coron and P.~Lissy.
\newblock Local null controllability of the three-dimensional Navier--Stokes
  system with a distributed control having two vanishing components.
\newblock {\em Inventiones mathematicae}, 198(3):833--880, 2014.


\bibitem{2018coronmarbach}
Jean--Michel  Coron, Fr{\'e}d{\'e}ric Marbach, and Franck Sueur.
\newblock Small--time global exact controllability of the Navier-Stokes equation with Navier 
	slip--with--friction boundary conditions
\newblock {\em arXiv preprint arXiv:1612.08087}, 2018.

\bibitem{FC-B-G-P}
E.~Fern{\'a}ndez-Cara, M.~Gonz{\'a}lez-Burgos, S.~Guerrero, and J.-P. Puel.
\newblock Null controllability of the heat equation with boundary Fourier
  conditions: the linear case.
\newblock {\em ESAIM: Control, Optimisation and Calculus of Variations},
  12(3):442--465, 2006.


\bibitem{fernandez2006some}
E.~Fern{\'a}ndez-Cara, S.~Guerrero, O.~Y. Imanuvilov, and J.-P. Puel.
\newblock Some controllability results for the n-dimensional Navier--Stokes and
  Boussinesq systems with n-1 scalar controls.
\newblock {\em SIAM journal on control and optimization}, 45(1):146--173, 2006.


\bibitem{fursikov1996controllability}
A.~V. Fursikov and O.~Y. Imanuvilov.
\newblock {\em Controllability of evolution equations}.
\newblock Number~34. Seoul National University, 1996.

\bibitem{fursikov1998local}
A.~V. Fursikov and O.~Y. Imanuvilov.
\newblock Local exact boundary controllability of the Boussinesq equation.
\newblock {\em SIAM Journal on Control and optimization}, 36(2):391--421, 1998.

\bibitem{fursikov1999exact}
A.~V. Fursikov and O.~Y. Imanuvilov.
\newblock Exact controllability of the Navier-Stokes and Boussinesq equations.
\newblock {\em Russian Mathematical Surveys}, 54(3):565--618, 1999.


\bibitem{guerrero2006local}
S.~Guerrero.
\newblock Local exact controllability to the trajectories of the Navier-Stokes
  system with nonlinear Navier-slip boundary conditions.
\newblock {\em ESAIM: Control, Optimisation and Calculus of Variations},
  12(3):484--544, 2006.

\bibitem{guerrero2006local_boussinesq}
S.~Guerrero.
\newblock Local exact controllability to the trajectories of the Boussinesq
  system.
\newblock In {\em Annales de l'IHP Analyse non lin{\'e}aire}, volume~23, pages
  29--61, 2006.
  

\bibitem{guerrero2018local}
S.~Guerrero and C.~Montoya.
\newblock Local null controllability of the n-dimensional Navier--Ntokes system
  with nonlinear Navier-slip boundary conditions and n- 1 scalar controls.
\newblock {\em Journal de Math{\'e}matiques Pures et Appliqu{\'e}es},
  113:37--69, 2018.



\end{thebibliography}
\end{document}